\documentclass[11pt,twoside]{article}
\usepackage{latexsym}
\usepackage{amssymb,amsbsy,amsmath,amsfonts,amssymb,amscd}
\usepackage{graphicx}
\usepackage{color}
\usepackage{hyperref}
\usepackage{caption}

\newcommand{\commentout}[1]{}
\newcommand{\R}{\mathbb{R}}
\newcommand{\N}{\mathbb{N}}

\newcommand {\e}  {\varepsilon}
\newcommand {\da} {\delta}

\newcommand {\vp} {\varphi}

\newcommand {\Chi} {{\bf \raise 2pt \hbox{$\chi$}} }
\newcommand {\f}   {\frac}
\newcommand {\p}   {\partial}

\newcommand {\proof} {\noindent {\bf Proof}. }
\newcommand{\beq}{\begin{equation}}
\newcommand{\eeq}{\end{equation}}
\newcommand{\bal}{\begin{align}}
\newcommand{\bc}{\begin{cases}}
\newcommand{\ec}{\end{cases}}
\newcommand{\bea} {\begin{array}{rl}}
\newcommand{\eea} {\end{array}}
\newcommand{\fer} {\eqref}
\newcommand{\bepa}{\left\{ \begin{array}{l}}
\newcommand{\eepa} {\end{array}\right.}
\newtheorem{theorem}{Theorem}[section]
\newtheorem{lemma}[theorem]{Lemma}

\newcommand{\qed}{{ \hfill
                     {\unskip\kern 6pt\penalty 500 \raise -2pt\hbox{\vrule\vbox to 6pt{\hrule width 6pt
                     \vfill\hrule}\vrule} \par}   }}
\title{\Large \bf A homogenization approach for the motion of motor proteins }

\author{Sepideh Mirrahimi  \thanks{
CMAP, Ecole Polytechnique, CNRS, INRIA. Route de Saclay, 91128 Palaiseau Cedex, France.  Email: mirrahimi@cmap.polytechnique.fr.}
\and Panagiotis E. Souganidis  \thanks{The University of Chicago,
Department of Mathematics, 5734 S. University Avenue,
Chicago, IL 60637, USA. Email: souganidis@math.uchicago.edu.} \thanks{Partially supported by the National Science
Foundation.}}

\date{\today}

\begin{document}
\maketitle
\pagestyle{plain}
\pagenumbering{arabic}

\begin{abstract}
We consider the asymptotic behavior of an evolving weakly coupled Fokker-Planck system of two equations set in a periodic environment. The magnitudes of the diffusion and the coupling are respectively proportional and inversely proportional to the size of the period. 
We prove that, as the period tends to zero, the solutions of the system either propagate (concentrate) with a fixed constant velocity (determined by the data) or do not move at all. The system arises in the modeling  of motor proteins which can take two different states. Our result implies that, in the limit,
the molecules either move along a filament with a fixed direction and constant speed or remain  immobile.
\end{abstract}


\noindent{\bf AMS Class. No:} 35B25, 35B27, 49L25, 92C05
%
%
\section{Introduction}
\label{intro}

We study the asymptotics, as $\e \to 0$, of the weakly coupled  Fokker-Planck system

\begin{equation}\label{eq:model}
\left\{ \begin{array}{ll}
n_{\e,t}^{1}-\e \Delta n_{\e}^{1}-\mathrm{div}_x(n_{\e}^{1}D_y \psi(\f{x}{\e}))+\f{1}{\e}\nu^{1}(\f{x}{\e})n_\e^{1}=\f{1}{\e}\nu^{2}(\f{x}{\e})n_\e^{2} \\
&\quad\text{in }\ \R^d\times(0,\infty),\\
 n_{\e,t}^{2}-\e \Delta n_{\e}^{2}+\f{1}{\e}\nu^{2}(\f{x}{\e})n_\e^{2}=\f{1}{\e}\nu^{1}(\f{x}{\e})n_\e^{1}.
\end{array}
\right.
\end{equation}

Systems like \eqref{eq:model} have been used to model the motion of motor proteins along molecular filaments or microtubules \cite{CD.BE.GO:95,FJ.AA.JP:97,RA.PH:02,DK.MK:02,MC.SH.DK:04,JD.DK.MK:04}.  The intracellular transport in eukaria is attributed to  motor proteins that transform the chemical energy into mechanical motion. For example myosins, which are known for their role in the muscle contraction, move along actin filaments and kynesins move along microtubules

In  \eqref{eq:model} the molecules have the conformations  $1$ and $2$ with densities $n_\e^{1}$ and $n_\e^{1}$ respectively and are
influenced by the periodic potential $\psi$ provided by the filaments (see \cite{AH:57}). The (periodic) functions $\nu^{1}$ and $\nu^{2}$ indicate the rates of change between the two states.
The existence of traveling waves, the asymptotic speed of propagation for large time as well as the presence of concentration effects for models with diffusion and a periodic drift have been studied in many papers; for instance see  \cite{PC.SM:08,AB.JD.MK:09}. A system similar to  \eqref{eq:model} also appears in the stochastic Stokes' drift where particles are suspended in a liquid and are subjected to diffusion and a net drift due to the presence of a wave in the liquid \cite{IB.CV.MC:00}. From the mathematical point of view, we refer to \cite{GA.YC:00} for a related homogenization problem as well as to  \cite{PD.AP:05} for the homogenization of a similar equation with parabolic scaling. 

To  formulate our result for the densities $n_\e^1$ and $n_\e^2$ next we introduce the assumptions we will be using throughout the paper. In addition to
\beq\label{ass:coef} \nu^{1}, \nu^{2} \ \text{and } \psi \ \text{ are smooth and } 1 \text{-periodic and } \  \nu^{1}>0,  \ \nu^{2}>0, \eeq
we assume that

\beq\label{as:pos}  
n_\e^1(x,0)>0
\ \text{ and } \  n_\e^2(x,0)>0  
 \ \text{ for all $x\in\R^d$ }, \eeq
\beq\label{mass:ini} I_\e^0=\int_\R n_\e^{1}(x,0)dx+\int_\R n_\e^{2}(x,0)dx \rightarrow I_0>0, \ \ \text{ as } \ \ \e\to0,\eeq
\beq\label{R2ini2}
\left\{ \begin{array}{ll}
\lim_{\e\rightarrow 0}\; \e\ln n_\e^1(\cdot,0)=\lim_{\e\rightarrow 0} \;\e \ln n_\e^2(\cdot,0)=-\infty \ \  \text{ locally uniformly in $\R^d\setminus\{0\}$, }\\[3mm]
\text{and}\\[2mm]
\underset{\underset{\e\rightarrow 0}{{y\rightarrow 0}}}\limsup\;\e\ln n_\e^i(y,0)\leq0 \ \text{ for } \ i=1,2,
\end{array}
\right.
\eeq
and, there exist constants  $A>0$ and  $B$ such that, for all $x\in\R^d$ and $i=1,2$,
\beq\label{as:massfini}
n_\e^{i}(x,0)\leq e^{\e^{-1}(-A|x|+B)} .
\eeq

Let $\delta$ denote the Dirac mass at the origin. Our first result is

\begin{theorem}\label{thm:trans}
Assume  \fer{ass:coef}, \fer{as:pos}, \fer{mass:ini}, \fer{R2ini2} and \eqref{as:massfini}. There exists $\bar v\in\R^d$ such that, as $\e\to0$ and in the sense of measures,
\beq\label{concentration}
n_\e^{1}(t,x)+n_\e^{2}(t,x)\rightharpoonup \delta(x-t\bar v)I_0. 
\eeq
\end{theorem}

To prove Theorem~\ref{thm:trans} we analyze the behavior, as $\e\to0$, of the functions $R_\e^1, R_\e^2:\R^d\times[0,\infty)\to\R$ which are obtained from $n_\e^1$ and $n_\e^2$ by the classical exponential change of variable (we show later in the paper that $n_\e^{1}>0$ and $n_\e^{2}>0$)

\beq\label{def:R}n_\e^{1}=\exp(-R_\e^{1}/\e) \ \ \text{ and } \ \   n_\e^{2}=\exp(-R_\e^{2}/\e),\eeq
and solve the system

\begin{equation}\label{eq:R}
\left\{ \begin{array}{ll}
R_{\e,t}^{1}-\e \Delta_x R_{\e}^{1}+|D_x R_{\e}^{1}|^2-D_y \psi(\f{x}{\e})\cdot D_x R_{\e}^{1}+\Delta_y \psi(\f{x}{\e})+\nu^{2}(\f{x}{\e})\exp(\f{R_\e^{1}-R_\e^{2}}{\e})=\nu^{1}(\f{x}{\e})\\
&\text{in }\ \R^d\times(0,\infty)\\
R_{\e,t}^{2}-\e \Delta_x R_{\e}^{2}+|D_x R_{\e}^{2}|^2+\nu^{1}(\f{x}{\e})\exp(\f{R_\e^{2}-R_\e^{1}}{\e})=\nu^{2}(\f{x}{\e}),
\end{array}
\right.
\end{equation}

where we use the subscript $x$ and $y$ to differentiate between differentiation with respect to $x$ and the fast variable $y=x/\e$.


To state the second main result we recall  that the ``half-relaxed'' upper and lower limits, denoted  $\overline R^i$ and $\underline R^{i}$ respectively, of the family $(R_\e^i)_{\e>0}$  are given, for $(x,t)\in\R^d\times[0,\infty)$ and $i=1,2$, by

\beq\label{liminfsup}\underline R^{i}(x,t)=\underset{\underset{\e\rightarrow 0}{{(y,s)\rightarrow (x,t)}}}\liminf R_\e^{1}(y,s) \ \ \text{and} \ \  \overline R^{i}(x,t)=\underset{\underset{\e\rightarrow 0}{(y,s)\rightarrow (x,t)}}\limsup R_\e^{1}(y,s).
\eeq


We have:

\begin{theorem}\label{thm:homo}
Assume \fer{ass:coef} and \fer{as:pos}. Then

(i) $\overline R^{1}=\overline R^{2}  \ \text{ and } \ \underline R^{1}=\underline R^{2} \ \text{in } \ \R^d\times (0,\infty).$

(ii) There exists a strictly convex $\overline H\in \mathrm{C}^1(\R^d)$ satisfying  $\overline H(0)=0$  and, for some $C>0$, $\overline H(p)\geq |p|^2-C$ such that $\overline R=\overline R^{1}=\overline R^{2}$ and $\underline R=\underline R^{1}=\underline R^{2}$ are respectively subsolution and supersolution of
\beq\label{effective}
\overline  R_t+\overline H( \overline R_x)\leq0 \ \ \text{ and } \ \  \underline  R_t+\overline H( \underline R_x)\geq0 \ \  \text{in } \R^d \times(0,\infty).
\eeq

(iii) Assume, in addition, \fer{R2ini2}. Then
\beq\label{lim:ini}
\overline R^{1}=\overline R^{2}=\underline R^{1}=\underline R^{2}=+\infty, \ \ \text{in }(\R^d\setminus\{0\}) \times\{0\},\quad  \underline R^{1}(0,0)\geq 0 \ \text{ and }  \  \underline R^2(0,0)\geq 0.\eeq

(iv) The special direction $\bar v$ in Theorem~\ref{thm:trans} is
\begin{equation}\label{direction}
\bar v =D \overline H(0).
\end{equation}
\end{theorem}

The asymptotic behavior of a time-independent version of \eqref{eq:model} set in $[0,1]$, which is also controlled by the same effective Hamiltonian
$\overline H$, was studied in  \cite{BP.PS:09} where it was proved  that, for asymmetric potentials,  the mass concentrates at either $x=0$ or $x=1$.  The asymmetry condition of \cite{BP.PS:09} is $D \overline H(0)\neq 0$. This behavior is in agreement  with our study of the time dependent problem. 
Indeed we prove here that, if $\bar v = D \overline H(0) \neq 0$,  the proteins (mass) move (spread) with constant velocity  $\bar v$. Hence the mass concentrates, for large times, on one end point of the filament, if the latter is assumed to have finite length. We  also refer to \cite{BP.PS:09bis} for a large deviation approach for the asymptotic behavior of the stationary solution of a similar model but with two potentials and to \cite{BP.PS:11} for the study of flashing ratchets.\\

By slight modifications of the proofs, all the results in this paper extend to systems  with two potentials like

\beq\label{2pot}
\left\{ \begin{array}{ll}
n_{\e,t}^{1}-\e \Delta n_{\e}^{1}-\mathrm{div}_x(n_{\e}^{1}D_y \psi^1(\f{x}{\e}))+\f{1}{\e}\nu^{1}(\f{x}{\e})n_\e^{1}=\f{1}{\e}\nu^{2}(\f{x}{\e})n_\e^{2} \\
&\quad\text{in }\ \R^d\times(0,\infty),\\
 n_{\e,t}^{2}-\e \Delta n_{\e}^{2}-\mathrm{div}_x(n_{\e}^{2}D_y \psi^2(\f{x}{\e}))+\f{1}{\e}\nu^{2}(\f{x}{\e})n_\e^{2}=\f{1}{\e}\nu^{1}(\f{x}{\e})n_\e^{1},
\end{array}
\right.
\eeq
the only difference being in the value of the effective Hamiltonian $\overline H(\cdot)$ --see Section \ref{cell}. We also remark that we cal also consider without any difficulty systems with more than two equations.\\

Our work is inspired from the ideas in wavefront propagation and large deviations \cite{CE.PS:89,GB.LE.PS:90}, the method of perturbed test functions in homogenization  \cite{LE:89} and the methods used in the study of the concentration effects \cite{GB.BP:07,GB.BP:08} and motor effects \cite{BP.PS:09bis, BP.PS:09,BP.PS:11}.\\

Throughout the paper solutions are taken to be either classical, if smooth or, otherwise, in the viscosity sense.
We refer to  \cite{C.I.L:92,GB:94} for a general introduction to the theory of the latter. In addition we denote by $C$ positive constants which are independent of $\e$ but may change from line to line. Moreover, $B_r(x)$ is the open ball in $\R^d$ centered at $x$ and of radius $r>0$ and $\overline B_r(x)$ stands for its closure. When $x=0$ we simply write $B_r$ and $\overline B_r$ respectively. Finally many statements in the paper hold for $i=1,2$ without any changes. Hence, unless necessary, we will not be repeating the ``for $i=1,2$''.

The paper is organized as follows. In Section \ref{cell} we introduce the cell problem corresponding to \fer{eq:R},  we recall that it has a solution and we introduce $\overline H$. In Section \ref{sec:boundR} we present some preliminary facts about the family $(n_\e^i)_{\e>0}$ and study the  properties of the families $(R_\e^{i})_{\e>0}$ that are needed to prove the convergence to $\overline R^i, \underline R^i$. Theorem~\ref{thm:homo} is proved in Section \ref{sec:conv}.
Using the results on the asymptotic behavior of the family $(R_\e^{i})_{\e>0}$ we prove Theorem~\ref{thm:trans} in Section~\ref{sec:concen}.  The asymptotic behavior of the family $(n_\e^{i})_{\e>0}$ in a more general setting is analyzed in Section \ref{sec:sev}. In Section \ref{sec:comp} we compare our results with the case of the parabolic scaling in \cite{PD.AP:05}. Finally, for the convenience of the reader, we present in the Appendix a sketch of the proof of the solvability of the cell problem and the properties of $\overline H$.

\section{The cell problem}
\label{cell}

In view of the presence of the exponential terms in \fer{eq:R} it is natural to expect that the $R_\e^1$'s and $R_\e^2$'s converge, if at all, as $\e\to0$ to the same limit $R$. Following \cite{BP.PS:09} we insert in \fer{eq:R} the formal expansion
$$R_\e^{i}(x,t)=R(x,t)+\e\phi^{i}(\f{x}{\e})+O^{i}(\e^2),$$
and, keeping only the terms multiplying $\e^0$ and writing $y$ for the fast variable $\f{x}{\e}$, we conclude that

\begin{equation}\label{eq:phi}
\left\{ \begin{array}{ll}
R_t-\Delta_y\phi^{1}+|D_y\phi^{1}+D_x  R|^2 -D_y\psi\cdot(D_y\phi^{1}+D_x R)+\Delta_y\psi+\nu^{2}\exp(\phi^{1}-\phi^{2})=\nu^{1}\\
&\text{in }\ \R^d\times(0,\infty)\\
R_t-\Delta_y\phi^{2}+|D_y\phi^{2}+D_xR|^2+\nu^{1}\exp(\phi^{2}-\phi^{1})=\nu^{2}.
\end{array}
\right.
\end{equation}

The goal is then to come up with  $\phi^{i}$'s so that \fer{eq:phi} is independent of $y$. This leads to the problem
to find, for each $p\in\R^d$, a unique constant $\overline H(p)$ such that the system, which is usually called the cell problem,

\begin{equation}\label{eq:cell}
\left\{ \begin{array}{ll}
-\Delta_y\phi^{1}+|D_y\phi^{1}+p|^2-D_y\psi\cdot(D_y\phi^{1}+p)+\Delta_y\psi+\nu^{2}\exp(\phi^{1}-\phi^{2})=\nu^{1}+\overline H(p) \\
&\text{in }\ \R^d, \\
-\Delta_y\phi^{2}+|D_y\phi^{2}+p|^2+\nu^{1}\exp(\phi^{2}-\phi^{1})=\nu^{2}+\overline H(p),
\end{array}
\right.
\end{equation}
admits an $1$-periodic solution $(\phi^{1}, \phi^{2})$ called the corrector.


We have:

\begin{lemma}\label{lem:cell} For each $p\in \R^d$ there exists a unique constant $\overline H(p)$ such that \fer{eq:cell}
has an $1$-periodic solution $(\phi^{1},\phi^{2})$. Moreover, $\overline H\in \mathrm{C}^1(\R^d)$, $\overline H(0)=0$, $\overline H$ is strictly convex, there exists a constant $C$ such that $\overline H(p)\geq |p|^2-C$, and , hence, $\overline H(p)\rightarrow \infty$ as $|p|\rightarrow \infty$.
\end{lemma}


A proof of Lemma~\ref{lem:cell} for $d=1$ was included in \cite{BP.PS:09}. For the convenience of the reader we sketch in the Appendix  the proof in $\R^d$.\\ 

When considering the system \fer{2pot} with two potentials, \fer{eq:cell} is replaced by
$$\left\{ \begin{array}{ll}
-\Delta_y\phi^{1}+|D_y\phi^{1}+p|^2-D_y\psi^1\cdot(D_y\phi^{1}+p)+\Delta_y\psi^1+\nu^{2}\exp(\phi^{1}-\phi^{2})=\nu^{1}+\overline H(p) \\
&\text{in }\ \R^d. \\
-\Delta_y\phi^{2}+|D_y\phi^{2}+p|^2-D_y\psi^2\cdot(D_y\phi^{2}+p)+\Delta_y\psi^2+\nu^{1}\exp(\phi^{2}-\phi^{1})=\nu^{2}+\overline H(p),
\end{array}
\right.$$
\section{Some preliminaries and the properties of $R_\e^{i}$}
\label{sec:boundR}

We summarize in the next lemma some basic properties of the families $(n_\e^{i})_{\e>0}$. They are the conservation of mass, the strict positivity of the $n_\e^{i}$'s, and a global upper bound yielding that, as $\e\to0$, there is very little mass at infinity.

We have:

\begin{lemma}\label{lem:mass}
(i) For all $t\geq0$ and $\e>0$,
 \beq\label{conservation}\int_{\R^d} n_\e^{1}(x,t)dx+\int_{\R^d} n_\e^{2}(x,t)dx=I_\e^0.\eeq
(ii) Assume \fer{ass:coef} and \fer{as:pos}. Then
\beq\label{npos}0<n_\e^i \ \text{ in } \ \R^d\times [0,+\infty).\eeq
(iii) Assume, in addition, \fer{as:massfini}. There exists $D>0$ such that, for all $(x,t)\in \R^d\times [0,+\infty)$.
 \beq\label{nborn} n_\e^i(x,t)\leq\exp(\f{-A|x|+B+Dt}{\e}).\eeq
In particular, there is small mass at infinity, i.e., for all $t\geq0$, there exists $M=M(t)>0$ such that
\beq \label{mass0}\int_{|x|\geq M} n_\e^{i}(t,x)dx \underset{\e\rightarrow 0}\longrightarrow 0.\eeq
\end{lemma}

\proof
The conservation of the mass follows from adding the equations in \fer{eq:model} and integrating over $\R^d$.

The form of \fer{eq:model} and \fer{ass:coef} and \fer{as:pos} allow us to use maximum principle-type arguments to obtain
\fer{npos} and \fer{nborn}.

Indeed let
$$F_{1,\e}(n)=n_{t}-\e \Delta n-\mathrm{div}_x(nD_y \psi(\f{x}{\e}))+\f{1}{\e}\nu^{1}(\f{x}{\e})n-\f{1}{\e}\nu^{2}(\f{x}{\e})n,$$
and
$$F_{2,\e}(n)=n_{t}-\e \Delta n+\f{1}{\e}\nu^{2}(\f{x}{\e})n-\f{1}{\e}\nu^{1}(\f{x}{\e})n.$$

It is easy to verify that $n_\e=\min(n_\e^1,n_\e^2)$ and $N_\e=\max(n_\e^1,n_\e^2)$ satisfy, in $\R^d\times(0,\infty)$, respectively

\beq\label{maxF} \max (F_{1,\e}(n_\e),{F_{2,\e}(n_\e)})\geq 0 \ \text{ and } \ \min (F_{1,\e}(N_\e),{F_{2,\e}(N_\e)})\leq 0.\eeq

Since $0$ is clearly a solution of the first inequality in \fer{maxF},  \fer{npos} follows from the strong maximum principle.

To prove \fer{nborn} we observe that, for $D$ sufficiently large,
$$G_\e(x,t)=\exp(\f{-A|x|+B+Dt}{\e})$$
is a viscosity supersolution of the second inequality in \fer{maxF}.
Since, in view of \fer{as:massfini}, we also have
$$N_\e(x,0)
\leq G_\e(x,0),$$
we conclude using again the comparison principle.


Finally \fer{mass0} follows from \fer{nborn} after an appropriate choice of the constant $M$.
\qed

\medskip

We turn now to the properties of the $R_\e^i$'s which are presented in Theorem~\ref{thm:bounds} below. The proof is rather long. The first part, which provides an one-sided Lipshitz-type continuity in time, is based on the classical Harnack inequality. The lower bound in part(ii) follows from part(i). The arguments leading to the upper bound (part(iii)) are more tedious and require as an intermediary step, namely, the construction, again using part(i), of an appropriate  local upper bound.

We have:
\begin{theorem}
\label{thm:bounds}
(i) Assume \fer{ass:coef} and \fer{as:pos}. For all $\da>0$, there exists $C_\da>0$, such that, for all $\e>0,\; |z-z'|\leq \e,\;  \e \da\leq t_0$ and $i,j=1,2$,
\beq\label{harnackR}R_\e^{j}(z',t_0+\e)-R_\e^{i}(z,t_0)\leq \e \,C_\da.
\eeq

(ii) Assume \fer{ass:coef}, \fer{as:pos}, \fer{mass:ini} and \fer{R2ini2}. 
For any $a\in(0,\infty)$, there exists  $\e_0=\e_0(a)>0$ such that, for all $\e\leq \e_0$,
\beq\label{h} R_\e^{i}\geq -a \ \text{ in } \ \R^d\times[0,\infty).\eeq

(iii) Assume \fer{ass:coef}, \fer{as:pos}, \fer{mass:ini} and \fer{as:massfini}. For any compact subset  $K$ of $\R^d\times (0,\infty)$, there exist  $C_K>0$ and $\e_1=\e_1(K)>0$ such that, for all $\e\leq \e_1$ and $(x,t)\in K$,
\beq\label{k} R_\e^{i}(x,t)\leq C_K.\eeq
\end{theorem}

\proof We begin with the  \newline
{\bf Proof of \fer{harnackR}:} Observe that 
$\widetilde{n}^{1}(y,\tau)={n}_\e^{1}(\e{y},\e{\tau}) \ \text{ and } \ \widetilde{n}^{2}(y,\tau)={n}_\e^{2}(\e{y},\e{\tau})$
are positive (recall \fer{npos}) solutions to 
\beq\label{eq:ntilde}
\left\{\begin{array}{ll}
\widetilde n_{\tau}^{1}- \Delta_y\widetilde n^{1}-\widetilde n^{1}\Delta_y\psi(y) -D_y\psi(y)\cdot D_y\widetilde  n^{1}+\nu^{1}\widetilde n^{1}=\nu^{2}\widetilde n^{2}\\
&\quad\text{in }\ \R^d\times(0,\infty)\\
\widetilde n_{\tau}^{2}-\Delta_y\widetilde n^{2}+\nu^{2}\widetilde n^{2}=\nu^{1}\widetilde n^{1},
\end{array}
\right.
\eeq
a linear parabolic system with bounded, according to \fer{ass:coef}, coefficients.  It follows
from the classical Harnack inequality \cite{JM.64} that, for each $\da>0$, there exists, an independent of $\e$, $C_\da>0$ such that
for all $y_0\in \R^d, \ \tau_0\geq \da$ and $i,j=1,2$,
\beq\label{harnackt} \sup_{z\in B_1(y_0)} \widetilde n^{i}(z,\tau_0)\leq C_\da \inf_{z\in B_1(y_0)} \widetilde n^{j}(z,\tau_0+1 ).\eeq

Rewriting \fer{harnackt} in terms of $n^1$ and $n^2$ and in the original variables $(x,t)$ we get 

\beq\label{harnackn} \sup_{z\in B_\e(x_0)}  n_\e^{i}(z,t_0)\leq C_\da\inf_{z\in B_\e(x_0)} n_\e^{j}(z,t_0+\e ),\qquad \text{for }(x_0,t_0)\in \R^d\times [\e\da,+\infty).\eeq

Finally using \fer{def:R} we obtain \fer{harnackR}.\\

We continue with the \newline
{\bf Proof of the uniform bounds from below:}
Arguing by contradiction we assume that for some $a>0$ there exist $\e_k\to0$ and $(x_k,t_k)\in\R^d\times[0,\infty)$ such that
$$\min \left(R_{\e_k}^1(x_k,t_k),R_{\e_k}^2(x_k,t_k)\right)<-a.$$
Newt observe that, in view of \fer{R2ini2} and \fer{as:massfini}, for $\e\leq \e_a$ with $\e_a$ small enough, we have
$$\min \left(R_{\e_k}^1(x_k,t_k),R_{\e_k}^2(x_k,t_k)\right)>-\f a 2.$$
Moreover $\min \left(R_\e^1,R_\e^2\right)$ is a supersolution to
$$\begin{array}{c}\max \left(R_{\e,t}^{1}-\e \Delta_x R_{\e}^{1}+|D_x R_{\e}^{1}|^2-D_y \psi(\f{x}{\e})\cdot D_x R_{\e}^{1}+\Delta_y \psi(\f{x}{\e})+\nu^{2}(\f{x}{\e})-\nu^{1}(\f{x}{\e}),\right.\\
\left. R_{\e,t}^{2}-\e \Delta_x R_{\e}^{2}+|D_x R_{\e}^{2}|^2+\nu^{1}(\f{x}{\e})-\nu^{2}(\f{x}{\e})\right) \geq 0,\end{array}$$
which admits $-\f a 2-ct$ as a subsolution provided $c$ is chosen sufficiently large.

It follows that, for $\e\leq \e_a$ and $t \geq 0$,
$$ \min \left(R_\e^1(\cdot,t),R_\e^2(\cdot,t)\right)\geq -\f a 2-ct, \qquad \text{in $\R^d$},$$
and therefore, if $\da=\f {a}{2c}$, for $\e\leq \e_a$,
$$ \min \left(R_\e^1,R_\e^2\right)\geq -a,  \qquad \text{in $\R^d\times [0,\da]$}.$$
As a result the sequence $(x_k,t_k)$ chosen at the beginning of the proof must satisfy $t_k > \da$.\\

Using \fer{harnackR} we deduce that there exists $C_1>0$ such that, as $k\to\infty$ and for all $x$ such that  $|x-x_k|\leq \e_k$ and $i=1,2$,
$$R_{\e_k}^{i}(x,t_k+\e_k)\leq R_{\e_k}^{i}(x_k,t_k)+ C_1\e_k\leq -a+ C_1\e_k. $$

It follows that
$$\int_{\R^d} n_{\e_k}^{i}(x,t_k+\e_k)dx\geq \int_{|x-x_k|\leq \e_k} e^{-\f{R_{\e_k}^{i}(x,t_k+\e_k)}{\e_k}}dx\geq |B_{\e_k}(x_k)|e^{\f{a}{\e_k}-C_1}.$$

The right hand side of this inequality blows up as $k\to\infty$, while the left hand side is bounded in view of \fer{conservation} and \fer{mass:ini}, again a contradiction. \\

The last part of the proof is devoted to the \newline
{\bf Proof of the uniform upper bounds on compact:} Fix a compact subset $K$ of $\R^d\times(0,\infty)$, observe that
$$t_0=\inf \{ s\in (0,\infty) : \text{ there exists } x\in \R^d \text{ such that } (x,s)\in K \}>0,$$
choose $t_1\in (0,t_0)$ and write $\bar t_1=t_1/2$ and $\bar t_2=t_1/4$.

It follows from \fer{mass:ini}, \fer{conservation} and \fer{mass0} that there exist $\e_1>0$ and $M>0$ both dependent on $t_1$ such that, if $\e\leq \e_1$, then 
$$\int_{|x|\leq M}\, n_\e^{1}(x,\bar t_1)+n_{\e}^2(x,\bar t_1)\, dx \geq \f{I_0}{2},$$
and, hence, in view of  \fer{def:R}, there exists some $a>0$ such that
$$ \min_{\underset{ i=1,2}{|x|\leq M}}R_\e^i(x,\bar t_1)\leq b=-\ln (\f{aI_0}{M^d}).$$

Assume next that the above minimum is achieved at some point $x_\e\in \bar B_{M}$ and for $i=i^*$.
Applying \fer{harnackR} $L=\left\lfloor \f{t_1}{2\e}\right\rfloor$ times 
with $i=i^*$, $j=2$, $x_0=x_\e$ and $\da=\bar t_1$, we find  some $C=C_\da>0$ such that, for all $x\in B_{L\e}(x_\e)$,
\beq\label{bound1} R_\e^2(x,\bar t_1 +L\e)\leq b+CL\e\leq  b+C\bar t_1.\eeq

Choose $\gamma \geq b+C\bar t_1$ and for some $\beta>0$ to be fixed below define $\phi^1:B_{\bar t_1}(x_\e) \times (0,\infty)\to \R$ by
$$\phi^1(x,t)=\f{1}{{\bar t_1}^2 -|x-x_\e|^2}+\beta t+\gamma$$

We claim that, for $\e\leq \e_2=\min(\bar t_2 ,\e_1)$,
\beq
\label{bound2}R_\e^2\leq \phi^1 \quad \text{in } \quad Q_\e=B_{\bar t_1}(x_\e)\times [\bar t_1 +L\e,+\infty),\eeq
and, therefore,
\beq\label{bound3}R_\e^2\leq \phi^1 \quad \text{in } \quad Q_\e^1=B_{\bar t_1}(x_\e)\times [t_1,+\infty).\eeq

To prove \fer{bound2} we first notice that, in view of \fer{bound1} and the choice of $\gamma$
$$R_\e^2(\cdot,\bar t_1 +L\e)\leq  b+C\bar t_1 \leq \gamma\leq \phi^1(\cdot,\bar t_1 +L\e) \quad \text{in} \quad B_{\bf t_1}(x_\e) $$.

Moreover, if $\beta$ is large enough, using \fer{ass:coef}, for all $(x,t)\in Q_\e$, we have
$$\phi_t^1-\Delta_x\phi^1+|D_x\phi|^2=\beta-\e\left(\f{2d}{({\bar t_1}^2-|x-x_\e|^2)^2}+\f{8|x-x_\e|^2}{({\bar  t_1}^2-|x-x_\e|^2)^3}\right)+\f{4|x-x_\e|^2}{({\bar t_1 }^2 -|x-x_\e|^2)^4}>D=\max_{y\in \R^d}\nu^2(y).$$

The inequality above and  \fer{eq:R} yield that, for all $(x,t)\in Q_\e$,
$$ \phi^1_t-\e\Delta_x \phi^1 +|D_x\phi^1|^2+\nu^{1}(\f{\cdot}{\e})\exp(\f{R_\e^{2}-R_\e^{1}}{\e})>\nu^{2}. $$

Since $R_\e^{2}\leq \phi^1$ on the parabolic boundary of $Q_\e$,  \fer{bound2} follows from the maximum principle.


For $\e\leq \e_2$ let $ Q^2=(\R^d\setminus B_{\bar t_2}(x_\e))\times (t_1,+\infty)$ and for positive constants $\alpha, \eta,\zeta$ to be fixed below consider the map $\phi^2: (\R^d\setminus B_{\bar t_2}(x_\e))\times (t_1,+\infty)\to\R$ given by
$$\phi^2(x,t)=\f{\alpha|x-x_\e|^2}{t-t_1}+\eta t+\zeta.$$

We claim that
\beq
\label{bound4}R_\e^2\leq \phi^2 \quad \text{in } \quad Q^2.\eeq

As above we will show that, if we choose $\alpha,\eta$ and $\zeta$ appropriately, $\phi^2$ is a supersolution in $Q^2$ of the equation satisfied by
$R_\e^2$ and is above  $R_\e^2$ on the parabolic boundary of $Q^2$.

To this end  notice that, in view of \fer{bound3}, we may select $\zeta$ and $\eta$ large enough so that, for  $\e\leq \e_2$,
$$R_\e^2\leq \eta t+\zeta\leq \phi^2  \quad \text{ on } \quad \p B_{\bar t_2}(x_\e)\times (t_1,\infty).$$

Moreover, for sufficiently large $\alpha,\eta$,  we have
$$\phi_t^1-\Delta_x\phi^1+|D_x\phi|^2=\eta-\f{\alpha|x-x_\e|^2}{(t-t_1)^2}-\f{2\e \alpha d }{t-t_1}+\f{4\alpha^2|x-x_\e|^2}{(t-t_1)^2}>D,$$
and, hence, in $Q^2$,
$$\phi^2_t-\e\Delta_x \phi^2+|D_x\phi^2|^2+\nu^{1}(\f{\cdot}{\e})\exp(\f{R_\e^{2}-R_\e^{1}}{\e})>\nu^{2}(\f{\cdot}{\e}).$$

Since clearly $R_\e^2(\cdot,t_1)\leq\phi^2(\cdot,t_1) \ \text{in} \ (\R^d\setminus B_{\bar t_2}(x_\e))$, using again the maximum principle  we find \fer{bound4}.


To conclude observe that \fer{bound3} and \fer{bound4} yield that the family $(R_\e^{2})_{\e>0}$ is uniformly bounded from above in any compact subset of $\R^d\times (t_1,\infty[$ and thus, in particular, on $K$, for $\e\leq \e_1$. Finally, using again \fer{harnackR}, we deduce that the family $(R_\e^{1})_{\e>0}$ is also uniformly bounded from above on $K$, for $\e\leq \min(\e_1,t_0-t_1)$.
\qed

\section{Convergence to the Homogenized equation-The proof of Theorem~\ref{thm:homo}}
\label{sec:conv}

Before we begin the proof, we remark that, in view of the claimed properties of $\overline H$ (convexity and coercivity),
the equation
\beq\label{eq:eff}
  R_t+\overline H(  R_x)=0 \ \ \text{in } \ \  \R^d \times(0,\infty),
\eeq
admits a comparison principle even for initial data taking the  value $\infty$ (see \cite{cls}). However, this and  Theorem \ref{thm:homo} do not lead to $\overline R\leq\underline R$. Indeed we show in the next section that, in addition to \fer{lim:ini} in Theorem \ref{thm:homo}, the conservation of mass yields
$$\underline R(0,0)=\underline R^1(0,0)=\underline R^2(0,0)=0.$$
To be able to use the comparison principle it is necessary to  have
$\overline R^{1}(0,0)=\overline R^{2}(0,0)\leq0$, which is not possible. Indeed, since $\overline R^1$ and $\overline R^2$ are upper semicontinuous, \fer{lim:ini} implies
 $\overline R^{1}(0,0)=\overline R^{2}(0,0)=+\infty$.

We continue with the

\noindent {\bf Proof of Theorem \ref{thm:homo}.} It is immediate from \fer{harnackR} that (i) holds. Moreover
it follows from Theorem \ref{thm:bounds} (ii) that
\beq\label{Rpos}0\leq \underline R^{1}\leq \overline R^{1}\ \text{ and } \ 0\leq \underline R^{2}\leq \overline R^{2} \ \text{on } \ \R^d\times [0,\infty),\eeq
and thus, in particular, the second part of \fer{lim:ini}.


To prove \fer{effective} we employ  the so called perturbed test function method. Since the arguments are similar here we only show that $\overline R$ is a subsolution of \fer{effective}.

To this end we assume that, for some smooth $\vp$,  $\overline R-\vp$ achieves  a strict local maximum at $(x_0,t_0)\in \R^d\times(0,\infty)$. Following \cite{BP.PS:09},  we perturb $\vp$ using the solution $(\phi^{1},\phi^{2})$ of the cell problem \fer{eq:cell} corresponding to $p=D \vp(x_0,t_0)$. It follows that there exists $(x_\e,t_\e)\in\R^d\times(0,\infty)\to(x_0,t_0)$, as $\e\to0$, such that the
$\max_{i=1,2} (R_\e^{i}-\vp-\e\phi^{i}(\f{\cdot}{\e}))$
is attained at $(x_\e,t_\e)$ and, without any loss of generality since the argument is identical, for $i=1$. Hence,
\beq\label{order} R_\e^{1}(x_\e,t_\e)-\e\phi^{1}(\f{x_\e}{\e})\geq R_\e^{2}(x_\e,t_\e)-\e\phi^{2}(\f{x_\e}{\e}).\eeq

That $R_\e^{1}$ is a solution of \fer{eq:R} yields 
$$\begin{array}{c}\vp_t(x_\e,t_\e)-\e\Delta_{x}\vp(x_\e,t_\e)-\Delta_y\phi^{1}(\f{x_\e}{\e})+|D_x\vp(x_\e,t_\e)+D_y\phi^{1}(\f{x_\e}{\e})|^2\\[2mm]
-D_y\psi(\f{x_\e}{\e})\cdot\left(D_x\vp(x_\e,t_\e)+D_x\phi^{1}(\f{x_\e}{\e})\right)
+\Delta_y\psi(\f{x_\e}{\e})+\nu^{2}(\f{x_\e}{\e})\exp (\f{R_\e^{1}-R_\e^{2}}{\e})\leq\nu^{1}(\f{x_\e}{\e}).
\end{array}$$

Using the latter and  \fer{order} and writing  $p_\e=D\vp(x_\e,t_\e)$ we get
$$\begin{array}{c}\vp_t(x_\e,t_\e)-\e\Delta_x\vp(x_\e,t_\e)-
\Delta_y\phi^{1}(\f{x_\e}{\e})+|p_{\e}+D_y\phi^{1}(\f{x_\e}{\e})|^2-D_y\psi(\f{x_\e}{\e})\cdot\left(p_{\e}+D_x\phi^{1}(\f{x_\e}{\e})\right)\\[2mm]
+\Delta_y\psi(\f{x_\e}{\e})+\nu^{2}(\f{x_\e}{\e})\exp (\phi^{1}(\f{x_\e}{\e})-\phi^{2}(\f{x_\e}{\e}))\leq\nu^{1}(\f{x_\e}{\e}).
\end{array}$$

It follows from the above inequality and the definition of the effective Hamiltonian \fer{eq:cell} that, for some $o(1)\to0$ as $\e\to0$,
$$\vp_t(x_\e,t_\e)-\e\Delta_x\vp(x_\e,t_\e)+\overline H(D_x\vp(x_\e,t_\e))\leq o(1),$$
and
the conclusion follows letting $\e\to0$.

%



To prove the  first part of \fer{lim:ini}, since the arguments are identical, here we only show that $\underline R^{1}(\cdot,0)=\overline R^{1}(\cdot,0)=\infty \ \text{ in } \ \R^d \setminus \{ 0 \}$. To this aim, we fix $\kappa>0$, we select an auxiliary function $\xi\in C^\infty(\R^d)$ such that
$
  0<\xi(x)<1 \ \text{for } \ x\in \R^d\setminus \{0\} \ \text{ and } \
  \xi(0)=0$,
and prove that
\beq\label{variational}\max(\, \underline R^1_t+\overline H(\underline R^1_x),\underline R^i-\kappa \xi\,)\geq 0 \qquad \text{in } \; \R^d\times [0,\infty).\eeq

 Since we already know that \fer{variational} holds in $\R^d\times(0,\infty)$, to conclude we assume that, for a smooth $\phi$,  $\underline R^1-\phi$ achieves a (strict) local maximum in $(x_0,0)$ and we prove that either
$$\underline R^1(x_0,0)\geq \kappa \xi(x_0),$$
or
\beq \label{super}\phi_t(x_0,0)+\overline H(\phi_x(x_0,0))\geq 0.\eeq

If $x_0=0$,  the former is clearly true since  $\underline R^1(0,0)\geq 0=\kappa\xi(0)$.

So we assume that $\underline R^1-\phi$ has a local maximum in $(x_0,0)$ with $x_0\neq 0$ and, in addition, that
$\underline R^1(x_0,0)<\kappa\xi(x_0)$. Repeating the arguments used earlier in the proof we obtain, for some $(x_\e,t_\e)\to(x_0,0)$ as $\e\to0$, we have
$$\vp_t(x_\e,t_\e)-\e\Delta_x\vp(x_\e,t_\e)+\overline H(D_x\vp(x_\e,t_\e))\geq o(1).$$
Indeed  \fer{R2ini2} and the facts that
$\underline R^1(x_0,0)<\kappa \xi(x_0)<\kappa$ and $\lim_{\e \rightarrow 0}R_\e^{1}(y,0)=+\infty$ for all $y$ near $x_0$ yield
that $t_\e>0$.  The claim now follows by letting $\e \to0$.


Assume next that, for some $x_0\neq 0$, $\underline R^1(x_0,0)=b<+\infty$. We fix $\da>0$ and let
$$\mu^\da(x,t)=-\f{|x-x_0|^2}{\da}-\gamma t$$
for $\gamma=\gamma(\da)>0$ to be determined later.

Since $\underline R^1$ is lower semicontinuous, $\underline R^1-\mu^\da$ attains a minimum at some $(x_\da,t_\da)\in \R^d\times [0,\infty)$ such that, as $\da\rightarrow 0$, $(x_\da,t_\da)\rightarrow (x_0,0)$, provided that $\gamma\to \infty$ as $\da\to0$.

Observe next that the choice of $(x_\da,t_\da)$ yields
\beq\label{Rdainf} \f{|x_\da-x_0|^2}{\da}\leq\underline R^1(x_\da,t_\da)+\f{|x_\da-x_0|^2}{\da}+\gamma t_\da\leq \underline R^1(x_0,0)=b,\eeq
and, hence,
$$|x_\da-x_0|\leq \sqrt{b\da}.$$

If $t_\da>0$, according to part (ii),  we must have
$\mu^\da_t(x_\da,t_\da)+\overline H(D \mu^\da(x_\da,t_\da))\geq 0$
and, hence,
\beq\label{contrad}-\gamma+\overline H\left(-\f{2(x_\da-x_0)}{\da}\right)\geq 0,\eeq
which cannot be true if we choose $\gamma>\sup_{|x|\leq \sqrt{b\da}}\overline H(-\f{2x}{\da})$.

Now we assume that $t_\da=0$. If
$\underline R^1(x_0,0)<\kappa \xi(x_0),$
then \fer{Rdainf} yields
$\underline R^1(x_\da,0)<\kappa_\da\xi(x_\da)$ for some $\kappa_\da\to\kappa$ as $\da\to0$.
Using \fer{variational} with $\kappa_\da$ at the point $(x_\da,0)$ we obtain again \fer{contrad} and thus a contradiction.

It follows that
$\underline R^1(x_0,0)>\kappa \xi(x_0),$
which also leads to a contradiction, since it holds for arbitrarily large $\kappa$ and $\underline R^1(x_0,0)=b<+\infty$.

The first part of \fer{lim:ini} now follows.
\qed
\section{The transport of the concentration points - The proof of Theorem~\ref{thm:trans}}
\label{sec:concen}

We present here the

\proof [Proof of Theorem \ref{thm:trans}]
Using \fer{effective}, \fer{lim:ini}, the standard optimal control formula \cite{PL:82,WF.HS:93,MB.ICD:96} and a barrier argument similar to the one of Section 5 in \cite{CE.PS:89}, we obtain that $\underline R$ satisfies
\beq \label{control}\underline R(x,t)\geq \inf_{\underset{\zeta(0)=0,\, \zeta(t)=x}{(\zeta(s),s)\in \R^d\times (0,\infty)+}}\int_0^t H^\star (\dot\zeta(s))ds+\max(\underline R^1(0,0),\underline R^2(0,0)),\eeq
with $H^\star(p)=\sup_{q\in \R^d}(p\cdot q-\overline H(q)).$

Observe next that, since  Lemma~$2.2$ yields
$$\lim_{|q|\rightarrow +\infty} p\cdot q-\overline H(q)\geq \lim_{|q|\rightarrow +\infty}p\cdot q-|q|^2+C=-\infty,$$
the maximum of $(p\cdot q-\overline H(q))$ is attained at some  $q_p\in\R^d$ such that
$p=D \overline H(q_p).$

Moreover, since $\overline H$ is strictly convex and $\overline H(0)=0$, we have
$$D\overline H(q_p)\cdot q_p-\overline H(q_p)\geq0 \ \text{with equality only if  } \ q_p=0.$$

Hence we deduce that $$H^*(p)>0 \ \text{ for all } \ p\neq D \overline H(0) \ \text{ and  } \ H^*\left(D \overline H(0)\right)=0.$$

Therefore, using \fer{control} and \fer{lim:ini}, we obtain that
\beq \label{tDH}\underline R\geq 0 \ \text{ in  } \ \R^d\times \R^+ \ \  \text{and } \ \ \text{if} \ \underline R(x,t)=0, \ \text{then} \  x=tD \overline H(0) \ \text{and } \ \underline R^1(0,0)=\underline R^2(0,0)=0.\eeq

It follows from the latter, \fer{def:R} and \fer{liminfsup} that, for $i=1,2$,
$$\underset{\underset{\e\rightarrow 0}{(y,s)\rightarrow (x,t)}}\limsup n_\e ^{i}(y,s) =0 \ \text{in } \ \R^d\times [0,+\infty)\setminus \{(t D \overline H(0),t)\,|\,t\in[0,+\infty)\}.$$

Finally the last claim and \fer{conservation} yield that the $(n_\e^{1}+n_\e^{2})$'s converge weakly along subsequences to a measure $n$  with
$$\mathrm{supp}\; n\subset\{(t\nabla \overline H(0),t)\,|\,t\in [0,\infty)\}.$$

Since we also know that, according to \fer{mass0}, no mass escapes to infinity as $\e\to0$,
we deduce \fer{concentration} using  \fer{conservation}. Moreover, in view of \fer{tDH}, we obtain that
$$\underline R^1(0,0)=\underline R^2(0,0)=0.$$
\qed

\section{The case with several Dirac masses initially}
\label{sec:sev}

We proved Theorem \ref{thm:trans} under assumptions \fer{R2ini2} which imply that the densities
$n_\e^{1}$ and $n_\e^{2}$ are both initially concentrated at the origin. The result can be generalized to
densities concentrated at several points and probably not on the same points. If this is the case, the initial condition
is written as

\beq\label{multidirac}
\left\{ \begin{array}{ll}
\lim_{\e\rightarrow 0}\; \e\ln n_\e^1(\cot,0)=\lim_{\e\rightarrow 0}\;  \e \ln n_\e^2(\cdot,0)=-\infty \ \  \text{locally uniformly in $\R^d\setminus\left(\mathcal{A}\cup \mathcal{B}\right)$, }\\[3mm]
\text{and}\\[2mm]
\underset{\underset{\e\rightarrow 0}{y\rightarrow x}}\limsup\;\e\ln n_\e^i(y,0)\leq0  \ \  \text{ for all $x\in \mathcal{A}\cup \mathcal{B}$}.
\end{array}
\right.
\eeq

We have:
\begin{theorem}\label{thm:multi}
Assume \fer{ass:coef}, \fer{as:pos}, \fer{mass:ini}, \fer{as:massfini} and \fer{multidirac}. Then, for $i=1,\,2$ and as $\e\to0$, along subsequences and in the sense of measures,
$
n_\e^{1}\rightharpoonup n^{1}$ and $ n_\e^{2}\rightharpoonup n^{2} with 
$
with
$$\mathrm{supp }\,(n^{1} +n^{2})(\cdot,t)\subset \mathcal{C}(t)=\{x_1+tD \overline H(0),\cdots,x_n+tD \overline H(0)\}\cup \{y_1+t D \overline H(0),\cdots, y_m+t D \overline H(0)\}.$$
\end{theorem}

\proof
The proof of Theorem \ref{thm:multi} follows along the same lines  as the one of Theorem \ref{thm:trans}. The only difference is that \fer{lim:ini} and \fer{control} are replaced respectively by
$$
\overline R^{1}=\overline R^{2}=\underline R^{1}=\underline R^{2}=+\infty \; \text{in }(\R^d\setminus \mathcal{C}(0)) \times\{0\}, \
\ \underline R^{1}\geq 0 \ \text{ and } \  \underline R^2\geq 0 \ \text{in }\mathcal{C}(0) \times\{0\},
$$
and
\beq \label{controlbis}\underline R(x,t)\geq \inf_{\underset{\zeta(0)\in \mathcal{C}(0),\, \zeta(t)=x}{(\zeta(s),s)\in \R^d\times (0,\infty)}}\left\{\int_0^t H^\star (\dot\zeta(s))ds+\max\left(\underline R^1(\zeta(0),0),\underline R^2(\zeta(0),0)\right)\right\}.\eeq

Therefore, we have
$$\underline R\geq 0 \ \text{ in } \ \R^d\times(0,\infty) \ \  \text{and } \ \  \text{ if } \ \ \underline R(x,t)=0, \ \text{then} \ \ x\in \mathcal{C}(t).$$

The other parts of the proof are similar.
\qed
%
%
%

\section{A comparison with results for parabolic scaling}
\label{sec:comp}
We describe here the connection between our result and the study in \cite{PD.AP:05} of the
asymptotics, as $\e\to0$, of the solutions to
\beq\label{ap}
n_{\e,t}=(a_{ij}(\f{x}{\e},\f{t}{\e^2})n_{\e,{x_j}})_{x_i}+\f{1}{\e} b_i(\f{x}{\e},\f{t}{\e^2})\,
n_{\e,{x_i}}+\f{1}{\e^2} c(\f{x}{\e},\f{t}{\e^2}) n_\e \ \text{ in } \ \R^d\times (0,\infty).\eeq

Notice that this is a single equation --- not a system---obtained after a parabolic scaling $(x,t)\mapsto (\f{x}{\e},\f{t}{\e^2})$ rather than the hyperbolic scaling $(x,t)\mapsto (\f{x}{\e},\f{t}{\e})$ of the problem we consider in this paper. Nevertheless, as we explain below, there are some similarities.

It is proved in \cite{PD.AP:05} that the solution $n_\e$ of \fer{ap} admits the expansion
\beq\label{piat}n_\e(x,t)=w(\f{x}{\e},\f{t}{\e^2})\exp(-\f{\lambda_0t}{\e^2})v^0(x-\f{\bar b}{\e}t,t)+o(1).\eeq

The result we obtain here with the hyperbolic scaling formally gives
\beq\label{compare}n_\e^{i}(x,t)=w^{i}(\f{x}{\e},\f{t}{\e})\rho^{i}(t)\da(x- D \overline H(0)t)+o(1).\eeq

In particular we obtain a Dirac mass instead of the function $v^0$ while the term $\exp(-\f{\lambda_0t}{\e^2})$ disappears because we have conservation of mass. Moreover $t/\e$ is replaced by  $t$, because of the difference in the scaling. There is, however, as we explain below  a close connection between the function $w$ in \fer{piat} and $w^{i}$ in \fer{compare}.

Indeed $w$ in \cite{PD.AP:05} is the principal eigenvector of the periodic cell problem
$$w_s-(a_{ij}(z,s)w_{z_j})_{z_i}-b_i(z,s)w_{z_i}-c(z,s)w=\Lambda_0w,$$
while we have
$$w^{i}=\exp(-{\phi^{i}}),$$
where $(\phi^1,\phi^2)$ is the principle eigenvector for the cell problem \fer{eq:cell} corresponding to $p=0$ with
the corresponding eigenvalue $\overline H(0)=0$. This is because for $(\bar x,\bar t)$ in the support of $n^{i}$, we have $R(\bar x,\bar t)=0$ and thus $D R(\bar x,\bar t)=0$.

\appendix

\section{The proof of Lemma~\ref{lem:cell}}
\label{app:cell}

{\bf Step $1$:} We prove that, for all $p\in \R^d$, there exists a unique $\overline H(p)$ with the properties stated in Lemma \ref{lem:cell}. Following \cite{BP.PS:09}, for $i=1,2$,
 we define
$$\Chi^{i}(y)=\exp(-p\cdot y-\phi^{i}(y)).$$ 

Multiplying the two equations in \fer{eq:cell} by $\Chi^{1}$ and $\Chi^{2}$ respectively, we obtain
 \beq
 \begin{cases}
-\Delta_y\Chi^{1}-\mathrm{div}_y(D_y\psi\; \Chi^{1})+\nu^{1}\Chi^{1}-\nu^{2}\Chi^{2}=-\overline H(p)\Chi^{1},\\[2mm]
-\Delta_y\Chi^{2}+\nu^{2}\Chi^{2}-\nu^{1}\Chi^{1}=-\overline H(p)\Chi^{2}.
\label{eq:Chi}
\end{cases}\eeq
with the boundary condition
\beq \label{period}y\to e^{p.y}\,\Chi^{i}(y) \ \text{is one-periodic and  } \ \Chi^{i}>0.\eeq

We also impose  the normalization
\beq\label{normal}\int_0^1 (\Chi^{1}(y)+\Chi^{2}(y) ) dy=1.\eeq

It follows from the Krein-Rutman Theorem that, for all $p\in \R^d$, there exists a unique constant $\overline H(p)$ and functions
 $(\Chi^{1}, \Chi^{2})$ satisfying \fer{eq:Chi} together with \fer{period} and \fer{normal}. \\

{\bf Step $2$:} To prove that $\overline H\in \mathrm{C}^1(\R^d)$ we rewrite \fer{eq:Chi} in terms of
$w^{i}_p(y)=e^{p\cdot y}\Chi^{i}(y)$, which in view of \fer{period} are $1$-periodic,  and obtain the new system

 \beq
 \begin{cases}\label{ss}
-\Delta_y w^{1}_p+2 p\cdot D_y w^{1}_p-\mathrm{div}_y(D_y\psi\; w^{1}_p)+\left(-|p|^2+p\cdot D_y  \psi+ \nu^{1}\right)w^{1}_p-\nu^{2}w^{2}_p=-\overline H(p)w^{1}_p,\\[2mm]
-\Delta_y w^{2}_p+2 p\cdot D_y w^{2}_p+\left(-|p|^2+\nu^{2}\right)w^{2}_p-\nu^{1}w^{1}_p=-\overline H(p)w^{2}_p.
\end{cases}\eeq

Assuming for the moment that $w^{1}_p$ and $w^{2}_p$ are differentiable with respect to $p$,  after differentiating \fer{ss} we get
$$
 \begin{cases}
-\Delta_y \p_p w^{1}_p+2 p\cdot D_y \p_p w^{1}_p-\mathrm{div}_y(D_y\psi\;\p_p w^{1}_p)+\left(-|p|^2+p\cdot \nabla_y  \psi+ \nu^{1}\right)\p_p w^{1}_p-\nu^{2}\p_p w^{2}_p\\
+2 D_y w^{1}_p+(-2p+D_y  \psi) w^{1}_p= -\overline H(p)\p_p w^{1}_p-\overline H'(p)w^{1}_p,\\[3mm]

-\Delta_y \p_p w^{2}_p+2 p\cdot D_y \p_p w^{2}_p+\left(-|p|^2+\nu^{2}\right)\p_p w^{2}_p-\nu^{1}\p_p w^{1}_p+2D_y  w^{2}_p -2pw^{2}_p=-\overline H(p)\p_p w^{2}_p\\-\overline H'(p) w^{2}_p.
\end{cases}
$$

Let $(w^{1}_{p,*},w^{2}_{p,*})$ be the solution to the adjoint system of  \fer{ss} --note that Fredholm's alternative implies the existence of such a solution. Multiplying the equations of \fer{ss} by $w^{1}_{p,*}$ and $w^{2}_{p,*}$, integrating with respect to $y$ and adding the two resulting equations we obtain
$$
{2\int w^{1}_{p,*}D_y w^{1}_pdy+2\int w^{2}_{p,*}D_y  w^{2}_pdy+\int D_y  \psi \,w^{1}_p w^{1}_{p,*}dy }
= (2p-D \overline H(p))(\int w^{1}_p w^{1}_{p,*} dy+\int w^{2}_p w^{2}_{p,*}dy),
$$
which yields a formula for $D \overline H(p)$. 

To prove the above claim rigorously, we write the system satisfied by difference quotients  with respect to $p$ and we use the same idea as above to prove that  the difference quotients $h^{-1}(\overline H(p+he_k)-\overline H(p))$ converge to the above formula, as $h\to0$. Indeed from this formulation we first obtain that $\overline H(p)$ is continuous with respect to $p$. Then we show that $w^{1}_p$ and $w^{2}_p$ are continuous with respect to $p$ and, finally, we pass to the limit $h\to0$.

To prove that $w^{1}_p$ and $w^{2}_p$  are continuous with respect to $p$, we show that, for any sequence $(p_n)_{n\in\N}$ such that $p_n\to p$, the  corresponding eigenfunctions $(w^{1}_{p_n},w^{2}_{p_n})$  converge to $(w^{1}_p,w^{2}_p)$  as $n\to\infty$.
This follows, for example, from the stability of viscosity solutions, the continuity of $\overline H$, the uniqueness (up to a multiplicative constant) of the eigenfunctions $(w^{1}_{p},w^{2}_{p})$ and the normalization condition \fer{normal}.\\

{\bf Step $3$:} The proof of $\overline H(0)=0$.

 The adjoint system to \fer{eq:Chi} is
\beq\label{sss}
\begin{cases}
 -\Delta_y u^{1}+D_y \psi \cdot D_y u^{1} + \nu^{1} u^{1}=\nu^{1} u^{2} -\overline H(p) u^{1},\\[2mm]
 -\Delta_y u^{2}+\nu^{2} u^{2}=\nu^{2} u^{1}-\overline H(p) u^{2},
\end{cases}
\eeq
with the condition
$$y\to e^{-p\cdot y}u^{i}(y) \ \text{$1$ -periodic and } \ u^{i}>0.$$

When $p=0$, \fer{sss} admits the trivial solution $(u^{1},u^{2})=(1,1)$ and the constant $\overline H(0)=0$.
The claim follows from the uniqueness of the positive eigenvector.\\

{\bf Step $4$:} The proof of the strict convexity  of $\overline H(p)$.

Arguing by contradiction we assume that there exist $p_1, p_2$ such
that
$$\f{\overline H(p_1)+\overline H(p_2)}{2}\leq \overline H(\f{p_1+p_2}{2}).$$

Let $(\phi_1^{1},\phi_1^{2})$ and $(\phi_2^{1},\phi_2^{2})$ be solutions of
\fer{eq:cell} corresponding to $p=p_1$ and $p=p_2$ respectively.  Adding the corresponding equations we find

$$\begin{array}{l}
-\Delta_y(\f{\phi_1^{1}+\phi_2^{1}}{2})+|D_y(\f{\phi_1^{1}+\phi_2^{1}}{2})+\f{p_1+p_2}{2}|^2-D_y\psi\cdot\left(D_y(\f{\phi_1^{1}+\phi_2^{1}}{2})+
\f {p_1+p_2}{2}\right)+\Delta_y\psi\\[2mm]
+\nu^{2}
\exp(\f{\phi_1^{1}+\phi_2^{1}-\phi_1^{2}-\phi_2^{2}}{2})\\[2mm]
<-\Delta_y(\f{\phi_1^{1}+\phi_2^{1}}{2})
+\f{1}{2}|D_y\phi_1^{1} +p_1|^2+\f{1}{2}|D_y\phi_2^{1} +p_2|^2
-D_y\psi\cdot\left(\f{D_y\phi_1^{1}+p_1}{2}\right)-D_y\psi\cdot\left(\f{D_y\phi_2^{1}+p_2}{2}\right)\\
+\Delta_y\psi+\f{\nu^{2}}{2}
\exp(\phi_1^{1}-\phi_1^{2})+\f{\nu^{2}}{2}
\exp(\phi_2^{1}-\phi_2^{2})\\[2mm]
=\nu^{1}+\f{\overline H(p_1)+\overline H(p_2)}{2})
\leq \nu^{1}+\overline H(\f{p_1+p_2}{2}),
\end{array}
$$
and
$$\begin{array}{l}
-\Delta_y(\f{\phi_1^{2}+\phi_2^{2}}{2})+|D_y(\f{\phi_1^{2}+\phi_2^{2}}{2})+\f{p_1+p_2}{2}|^2
+\nu^{1}
\exp(\f{\phi_1^{2}+\phi_2^{2}-\phi_1^{1}-\phi_2^{1}}{2})<\nu^{2}+\overline H(\f{p_1+p_2}{2}).\\
\end{array}$$

It follows that the pair $(\f{\phi_1^{1}+\phi_2^{1}}{2},\f{\phi_1^{2}+\phi_2^{2}}{2})$ is a strict subsolution to the cell problem
\fer{eq:cell} corresponding to $p=\f{p_1+p_2}{2}$. This, however, contradicts the fact that $\overline H(\f{p_1+p_2}{2})$ is
the principal eigenvalue of the system.\\

{\bf Step $5$:} The proof of $\overline H(p) \geq |p|^2 -C$.

Rewriting \fer{eq:cell} we find
$$
 \begin{cases}
-\Delta_y\phi^{1}+|D_y\phi^{1}|^2-(D_y\psi-2p)\cdot D_y\phi^{1}-D_y\psi\cdot p+\Delta_y\psi+\nu^{2}\exp(\phi^{1}-\phi^{2})=\nu^{1}+\overline H(p)-|p|^2,\\
-\Delta_y\phi^{2}+|D_y\phi^{2}|^2+2p\cdot D_y\phi^{2}+\nu^{1}\exp(\phi^{2}-\phi^{1})=\nu^{2}+\overline H(p)-|p|^2.
 \end{cases}$$

Assume next that $\max_{[0,1]}(\phi^{1},\phi^{2})$ is attained at the point $\bar y$. If
 $\phi^{1}(\bar y)>\phi^{2}(\bar y)$, then
$$-\Delta \phi^{1}(\bar y)\geq 0 \ \text{ and } \ D \phi^{1}(\bar y)=0,$$
and, hence,
$$\overline H(p) \geq -\|\Delta \psi\|_{L^{\infty}}-|p|\,\|D_y \psi\|_{L^{\infty}}-\|\nu^{1}\|_{L^\infty}+|p|^2.$$

If  $\phi^{1}(\bar y)>\phi^{2}(\bar y)$, then, similarly,  we find
$$\overline H(p) \geq -\|\nu^{1}\|_{L^\infty}+|p|^2.$$

The claim now follows.

\bigskip

\bigskip
\noindent {\bf Acknowledgment} This work has partly been done during a visit of the first author to the University of Chicago while supported by {\em Fondation Sciences Math\'ematique de Paris}.  S.~M. would like to thank  the University of Chicago for its hospitality.



\end{document}